\documentclass[preprint,12pt]{elsarticle}

\newtheorem{theorem}{\textbf{Theorem}}

\newtheorem{lemma}{\textbf{Lemma}}

\usepackage{mathrsfs, amsmath, amssymb}
\usepackage{mathtools}

\def\a {\alpha}
\def\b {\beta}

\def\N {\mathbb{N}}

\def\Q {\mathbb{Q}}
\def\R {\mathbb{R}}

\def\e {\epsilon}

\def\g {\gamma}

\usepackage{amssymb}

\journal{-}

\begin{document}

\begin{frontmatter}



\title{On an upper bound of the degree of polynomial identities regarding linear recurrence sequences}

\author[ap]{Ana Paula Chaves\fnref{app}\corref{AP}}{\ead{apchaves@ufg.br}}
\author[gg]{Carlos Gustavo Moreira}{\ead{gugu@impa.br}}
\author[en]{Eduardo Henrique Rodrigues do Nascimento}
\address[ap]{Instituto de Matem\' atica e Estat\'istica, Universidade  Federal de Goi\'as, Goi\^ania, 74001-970, Brazil}
\address[gg]{Instituto de Matem\' atica Pura e Aplicada, Rio de Janeiro, 22460-320, Brazil}
\address[en]{Arena Highschool, Goi\^ania, 74215-160, Brazil}

\fntext[app]{Supported in part by FAPERJ-Brazil}
\cortext[AP]{Corresponding author}

\begin{abstract}
Let $(F_n)_{n\geq 0}$ be the Fibonacci sequence given by $F_{n+2}=F_{n+1}+F_n$, for $n\geq 0$, where $F_0=0$ and $F_1=1$. There are several interesting identities involving this sequence such as $F_n^2+F_{n+1}^2=F_{2n+1}$, for all $n\geq 0$. Inspired by this naive identity, in 2012, Chaves, Marques and Togb\' e proved that if $(G_m)_m$ is a linear recurrence sequence (under weak assumptions) and $G_n^s+\cdots +G_{n+k}^s\in (G_m)_m$, for infinitely many positive integers $n$, then $s$ is bounded by an effectively computable constant depending only on $k$ and the parameters of $(G_m)_m$. In this paper, we  generalize this result, proving, in particular, that if $(G_m)_m$ and $ (H_m)_m$ are li\-near recurrence sequences (also under weak assumptions), $R(z) \in \mathbb{C}[z]$, and $ R(G_n)+\cdots +R(G_{n+k})$ belongs to $(H_m)_m$, for infinitely many positive integers $n$, then the degree of $R(z)$ is bounded by an effectively computable constant depending only on the parameters of $(G_m)_m$ and $(H_m)_m$ (but surprisingly not on $k$).

\end{abstract}


\end{frontmatter}



\section{Introduction}

A sequence $(G_n)_{n\geq 0}$ is said to be a \textit{linear recurrence sequence} with coefficients $c_0$, $c_1, \ldots, c_{k-1}$, where $c_0\neq 0,$ if
\begin{equation}\label{recu}
G_{n+k}=c_{k-1}G_{n+k-1}+\cdots + c_1G_{n+1}+c_0G_n,
\end{equation}
for all $n \geq 0$. A linear recurrence sequence is therefore completely determined by its \textit{initial values} $G_0,\ldots ,G_{k-1}$, and by the coefficients $c_0,c_1,\ldots,c_{k-1}$. The integer $k$ is called the {\it order} of the recurrence. The \textit{characteristic polynomial} of $(G_n)_{n\geq 0}$ is given by
$$G(x)=x^{k}-c_{k-1}x^{k-1}-\cdots - c_1x-c_0.$$
The roots of $G(x)$, named $r_1,r_2,\ldots ,r_{\ell}$, are called the \emph{roots of the recurrence}. A root $r_j$ of the recurrence is said a {\it dominant root} if $|r_j|>|r_i|$, for all $j\neq i\in \{1,...,\ell\}$. In this paper, we consider only recurrence sequences whose coefficients and initial values are real algebraic numbers, i.e., whose terms are real algebraic numbers. 

A general Lucas sequence $(C_n)_{n\geq 0}$ given by $C_{n+2}=aC_{n+1}+bC_n$, for $n\geq 0$, where the values $a,\ b,\ C_0$ and $C_1$ are previously fixed, is an example of a linear recurrence of order $2$ (also called {\it binary}). For instance, if $C_0=0$ and $C_1=a=b=1$, then $(C_n)_{n\geq 0}$ becomes the well-known \textit{Fibonacci sequence} $(F_n)_{n\geq 0}$, whose first terms are:
$$
0,1,1,2,3,5,8,13,21,34,55,89,144,\ldots .
$$

Fibonacci numbers are famous not only for featuring  in many unexpected branches of Mathematics, Computer Science,  Architecture,  Biology, and so on, but also (and mainly) for their stunning properties (see \cite{pos} and \cite{kalman}). Among many identities satisfied by these numbers, we have the following
\begin{equation}\label{naive}
F_n^2+F_{n+1}^2=F_{2n+1},\ \mbox{for\ all}\ n\geq 0,
\end{equation}
which tells us that the sum of the square of two consecutive Fibonacci numbers is still a Fibonacci number. Marques and Togb\'{e} \cite{MT} questioned what about such sums with higher powers, and showed that, if $x \geq 1$ is an integer such that $F_{n}^x + F_{n+1}^x$ is a Fibonacci number for all sufficiently large $n$, then $x \in \{1, 2 \}.$ In 2011, Luca and Oyono \cite{luca} solved this problem completely, proving that the  Diophantine equation
\begin{equation} \label{eq1}
F_m^s+F_{m+1}^s=F_n ,
\end{equation}
has no solutions $(m,n,s)$, for $m\geq 2$ and $s \geq 3$. Since then, many authors have considered and worked on variations of \eqref{eq1}, looking for solutions when the Fibonacci sequence is replaced by other sequences that, in a way, are high-order generalizations of $(F_n)_{n\geq 0}$, such as the \emph{$k$-generalized Fibonacci sequence} (see \cite{gersica, chaves, ruiz}). 

Our main interest relies on the following result, yet related to \eqref{naive}. In 2010, Chaves, Marques and Togb\'e \cite{CMT} considered the sum of many powers of terms of a linear recurrence sequence, and their main theorem stated that if $(G_n)_{n}$ is an integer linear recurrence sequence (under some assumptions), $s, k$ and $b$ are positive integer numbers and  $\epsilon_j \in \{0,1\}$, with $0 \leq j \leq k-1$, such that
\begin{equation}\label{rec1}
G_n^s + \epsilon_{1} G_{n+1}^s + \cdots + \epsilon_{k-1}G_{n+k-1}^s + G_{n+k}^s
\end{equation}
belongs to the sequence $(b \cdot G_n)_n$ for infinitely many $n \in \mathbb{N}$, then $s$ is bounded by an effectively computable constant depending on $k$, $b$ and the parameters of $(G_n)_n$. The assumption required for $(G_n)_n$ is that its characteristic polynomial has a simple positive root being the unique zero outside the unit circle. For instance, $(F_n)_n$ and some of its high-order generalizations satisfy this condition.  

The aim of this work is to extend the main result of \cite{CMT}, reaching a much more general set of linear recurrence sequences and also combinations of powers of terms of such sequences that could belong to other recurrence sequences. More precisely, we prove the following:

\begin{theorem}\label{main}
Let $(G_n)_{n \geq 0}$ and $(H_n)_{n \geq 0}$ be unbounded linear recurrence sequences whose terms are real algebraic numbers, having simple real dominant roots $\a$ and $\b$, respectively, such that $\a = \gamma^a$ and $\b = \gamma^b$, for $a,b \in \N$. Also, let $c \in \R \cap \overline\Q$, and $M,m \in \mathbb{R}_{>0}$. If $R(z)\in \mathbb C[z]$, is a  polynomial of degree $s$ whose leading coefficient is a real algebraic number, $k\in \N$, and for $0\leq j\leq k-1$ we have $\e_j \in \mathbb C$, not all zero, such that $m\le |\e_j|\le M$ for all $\e_j \neq 0$, and
$$ \e_0R(G_n)+\e_1R(G_{n+1})+\cdots +\e_{k-1}R(G_{n+k-1})+R(G_{n+k})$$
belongs to the sequence $(c \cdot H_n)_{n}$, for infinitely many $n$'s, then there exists an effectively computable constant $E$, such that $s<E$. This constant $E$ depends only on $c, m, M$ and the parameters of $(G_n)_n$ and $(H_n)_n$ (but not on $k$).
\end{theorem}

\section{Auxiliary results}\label{sec2}

In this section, we recall some results that will be fundamental for the proof of Theorem \ref{main}. Let $G(x)$ be the characteristic polynomial of a linear recurrence $(G_n)_n$. One can factor $G(x)$ over the set of complex numbers as
$$
G(x)=(x-r_1)^{m_1}(x-r_2)^{m_2}\cdots (x-r_{\ell})^{m_{\ell}},
$$
where $r_1,...,r_{\ell}$ are the roots of the recurrence and $m_1, \ldots ,m_{\ell}$ are positive integers.   A fundamental result in the theory of recurrence sequences asserts that there exist uniquely determined non-zero polynomials $g_1,...,g_{\ell}\in \mathbb{Q}(\{r_j\}_{j=1}^{\ell})[x]$, with $\deg g_j\leq m_j-1$, for $j=1,...,\ell$, such that
\begin{equation}\label{written}
G_n=g_1(n)r_1^n+\cdots +g_{\ell}(n)r_{\ell}^n,\ \mbox{for\ all}\ n.
\end{equation}
For more details, see \cite[Theorem C.1]{shorey}. We recall that, since $G(x) \in \overline{\mathbb{Q}}[x]$, then $g_j(n)$ is an algebraic number, for all $j=1,...,{\ell}$, and $n\in \mathbb{N}$. In the case of the Fibonacci sequence, \eqref{written} is known as {\it Binet's formula}:
\begin{center}
$F_n=\displaystyle\frac{\alpha^n-\beta^n}{\alpha-\beta}$,
\end{center}
where $\alpha=(1+\sqrt{5})/2$ (the golden number) and $\beta=(1-\sqrt{5})/2=-1/\alpha$.
If $r_j$ is a  dominant root  of the recurrence the corresponding polynomial $g_j(n)$ is named the \textit{dominant polynomial}. 

The formula given in \eqref{written} and some calculations will allow us to obtain a linear form in three logarithms and then determine lower bounds {\it \`a la Baker} for it. For that purpose, we deduce from the main result of Matveev \cite{matveev}, the following lemma.

\begin{lemma}\label{lemma1}
Let $\alpha_1,\alpha_2,\alpha_3$ be real algebraic numbers and let $b_1,b_2,b_3$ be non-zero integer rational numbers. Define
$$\Gamma=\alpha_1^{b_1}\alpha_2^{b_2}\alpha_3^{b_3}-1.$$
Let $D$ be the degree of the number field $\Q(\alpha_1,\alpha_2,\alpha_3)$ over $\Q$ and let $A_1,A_2,A_3$ be positive real numbers which satisfy
\begin{center}
$A_j\geq \max\{Dh(\alpha_j),|\log \alpha_j|,0.16\}$, for $j=1,2,3$.
\end{center}
Assume that $B\geq \max\{|b_j|; 1\leq j\leq 3\}$. Define also
$$
C_1=1.4\times30^6\times3^{4.5}\times D^2\log (eD)
$$
If $\Gamma\neq 0$, then
$$
|\Gamma|> \exp(-C_1A_1A_2A_3\log(eB)).
$$
\end{lemma}
As usual, in the previous statement, the \textit{logarithmic height} of an $n$-degree algebraic number $\tau$ is defined as $$h(\tau)=\frac{1}{n}(\log |a|+\sum_{j=1}^n\log \max\{1,|\tau^{(j)}|\}),$$
where $a$ is the leading coefficient of the minimal polynomial of $\tau$ (over $\mathbb{Z}$) and $(\tau^{(j)})_{1\leq j\leq n}$ are the conjugates of $\tau$.

The next lemma plays an important role in the first steps of the proof of Theorem \ref{main}. It can be deduced from \cite[Lemma 2]{CMT}, but we shall prove it quickly for the reader's convenience.

\begin{lemma}\label{lemma2}
Let $(G_n)_{n}$ be a linear recurrence such that its characteristic polynomial has a simple dominant root $r_1$. Then the dominant polynomial of $(G_n)_{n}$ is a non-zero real number $g_1$, and $G_n \sim g_1 r_1^n$. 
\end{lemma}

\noindent
\textbf{Proof.}
Since $r_1$ is a simple dominant root, it follows that $m_1-1=0$, thus $g_1(x)=g_1$ is constant. By dividing both sides of \eqref{written} by $r_1^n$, we note that 

\[\frac{G_n}{r_1^n} = g_1+\sum_{j = 2}^{\ell} g_j(n).\left(\frac{r_j}{r_1}\right)^n.\]
However, since $|r_j/r_1|<1$ for all $2 \leq j \leq \ell$, we have 
\[
\lim_{n \to \infty}g_j(n).\left(\frac{r_j}{r_1}\right)^n=0,
\]
therefore, having in mind that $r_1$ and $G_n$ are real numbers, we are done.
\qed

The following fact will be highlighted but stated without proof. 
\begin{lemma} \label{lemma3}
Let $r\in \mathbb R$ and $(a_n)_n$ be a sequence of integers. If
\[\lim_{n\to \infty} r^{a_n} = L,\]
then $L$ is either equal to zero or to an integer power of $r$.
\end{lemma}


Now we have all the ingredients needed to start dealing with the proof of Theorem \ref{main}, which follows in the next section. 

\section{The proof of Theorem \ref{main}} \label{sec3}

Throughout this section, we will use Landau's big-oh symbol, having in mind that the parameters  $c, m, M, (G_n)_n,(H_n)_n$ are given (meaning they are fixed), while $k,s$ are variables. 

First, we point out that both dominant roots $\a$ and $\b$ are outside the unit circle. Indeed, let $g$ and $h$ be the dominant polynomials of $(G_n)_n$ and $(H_n)_n$ respectively. Since from Lemma \ref{lemma2} we have $G_n \sim g\a^n$, if $|\a|\le 1$, then $(G_n)_n$ would be bounded, thus $|\a|>1$ and the same argument gives $|\b| >1$.

Now suppose that there exists an infinite subset $\mathcal{N} \subseteq \mathbb{N}$ and a sequence of positive integers $(t_n)_{n\in \mathcal N}$, such that for every $n$ in $\mathcal N$, we have 
\begin{equation} \label{maineq}
\e_0R(G_n)+\e_1R(G_{n+1})+\cdots+\e_{k-1}R(G_{n+k-1})+R(G_{n+k})= cH_{t_n} \ .
\end{equation}
Note that we can consider $R(z)$ to be monic, since the leading coefficient is a real algebraic number which can be ``embedded" into $c$, by dividing both sides of \eqref{maineq} by it.
Afterwards, we divide \eqref{maineq} by $\a^{sn}$, to get
\begin{equation}
\label{eqn:orig}
 \frac{ \e_0R(G_n)}{\a^{sn}} +\cdots+ \frac{\e_{k-1}R(G_{n+k-1})}{\a^{sn}}+\frac{R(G_{n+k})}{\a^{sn}}=\frac{cH_{t_n}}{\a^{sn}} \ .   
\end{equation}
Now, let us examine each part of the sum on the left-hand side of \eqref{eqn:orig}, as $n\in \mathcal N$ goes to infinity, and recall that $R(z)$ is monic and of degree $s$, to obtain for $i \in \{0,1, \ldots, k\}$,
\[\lim_{\begin{smallmatrix} n \to \infty & \\ n \in \mathcal{N} \end{smallmatrix}}\frac{R(G_{n+i})}{\a^{sn}}=\left[\lim_{\begin{smallmatrix} n \to \infty & \\ n \in \mathcal{N} \end{smallmatrix}} \frac{R(G_{n+i})}{G_{n+i}^s}\right] \cdot 
\left[\lim_{\begin{smallmatrix} n \to \infty & \\ n \in \mathcal{N} \end{smallmatrix}}\frac{G_{n+i}^s}{\a^{sn}}\right]=\lim_{\begin{smallmatrix} n \to \infty & \\ n \in \mathcal{N} \end{smallmatrix}}\left(\frac{G_{n+i}}{\a^{n}}\right)^s .\]
Since $\alpha$ is the dominant root of $(G_n)_n$, we use Lemma \ref{lemma2} once more to get
\[\lim_{\begin{smallmatrix} n \to \infty & \\ n \in \mathcal{N} \end{smallmatrix}}\left(\frac{G_{n+i}}{\a^{n}}\right)^s=\left(\lim_{\begin{smallmatrix} n \to \infty & \\ n \in \mathcal{N} \end{smallmatrix}}\frac{G_{n+i}}{\a^{n}}\right)^s=\a^{si}\cdot \left(\lim_{\begin{smallmatrix} n \to \infty & \\ n \in \mathcal{N} \end{smallmatrix}}\frac{G_{n+i}}{\a^{n+i}}\right)^s= g^{s}\a^{si},\]
and then the limit of the left-hand side of \eqref{eqn:orig} exists and equals to
\[g^s(\e_0+\e_1\a^s+\cdots+\e_{k-1}\a^{(k-1)s}+\a^{ks}) .\]
Therefore, the limit of the right-hand side of \eqref{eqn:orig} given by
\[
\lim_{\begin{smallmatrix} n \to \infty & \\ n \in \mathcal{N} \end{smallmatrix}} \frac{cH_{t_n}}{\a^{sn}} ,
\]
exists. Now, since $\a = \g^a$, we rewrite this limit and make use of $\b = \g^b$ and Lemma \ref{lemma2} to get,
\[
\lim_{\begin{smallmatrix} n \to \infty & \\ n \in \mathcal{N} \end{smallmatrix}} \frac{cH_{t_n}}{\g^{asn}} = \left(\lim_{\begin{smallmatrix} n \to \infty & \\ n \in \mathcal{N} \end{smallmatrix}} \frac{cH_{t_n}}{\b^{t_n}}\right) \cdot \left(\lim_{\begin{smallmatrix} n \to \infty & \\ n \in \mathcal{N} \end{smallmatrix}}\frac{\g^{bt_n}}{\g^{asn}}\right)=c h \cdot \lim_{\begin{smallmatrix} n \to \infty & \\ n \in \mathcal{N} \end{smallmatrix}}\g^{bt_n-asn} .
\]

Since it exists, Lemma \ref{lemma3} implies that we have two cases for the limit on the right-hand side. 
    
If $\displaystyle \lim_{\begin{smallmatrix} n \to \infty & \\ n \in \mathcal{N} \end{smallmatrix}}\g^{bt_n-asn} = 0$, then  
\[g^s(\e_0+\e_1\a^s+\cdots+\e_{k-1}\a^{(k-1)s}+\a^{ks}) = c h \cdot \lim_{\begin{smallmatrix} n \to \infty & \\ n \in \mathcal{N} \end{smallmatrix}}\g^{bt_n-asn} = 0 
 .\]
Thus, 
\[
\e_0+\e_1\a^s+\cdots+\e_{k-1}\a^{(k-1)s}+\a^{ks} = 0 ,
\]
hence by the triangle inequality and the fact that $|\e_j|\le M$, for $0 \leq j \leq k-1$, we get
\[ |\a|^{ks} \le |\e_0|+|\e_1||\a|^s+\cdots+|\e_{k-1}||\a|^{(k-1)s}  < M \frac{|\a|^{ks}}{|\a|^s-1} ,\]
which gives
\[s < \frac{\log (M+1)}{\log |\a|},\]
so our constant can be taken as $E:= \log (M+1)/\log |\a|$, and we are done. 

For the other case, we have
\begin{equation} \label{case2}
g^s(\e_0+\e_1\g^{as}+\cdots+\e_{k-1}\g^{(k-1)as}+\g^{kas}) = c h\g^{t}.
\end{equation}
where $t$ in an integer number given by $t=\lim_{n\to \infty}(bt_n-asn)$. Before going any further, we prove the following claim, which will be crucial to deduce that $s$ is bounded in our final step.

{\bf Claim:} \emph{$|t-aks| = O(s)$.}

Indeed, by dividing both sides of \eqref{case2} by $g^s\g^{kas}$, we obtain
\begin{equation} \label{case2.2}
    \e_0\g^{-kas}+\e_1\g^{-(k-1)as}+\cdots+\e_{k-1}\g^{-as}+1 = chg^{-s}\g^{t-kas} \ .
\end{equation}
Take $\Gamma :=  \e_0\g^{-kas}+\e_1\g^{-(k-1)as}+\cdots+\e_{k-1}\g^{-as}$. First, we have the following upper bound for $|\Gamma|$, which will be used a few times from now on:
\[
|\Gamma|\leq |\e_0||\g|^{-kas}+\cdots+|\e_{k-1}||\g|^{-as} \le M \cdot \frac{1-|\g|^{-kas}}{|\g|^{as}-1} ,
\]
\begin{equation} \label{desgamma}
    \Rightarrow \ |\Gamma| < \frac{M}{|\g|^{as}-1} \ .
\end{equation}
From now on, we may assume that $|\Gamma|<1/2$, otherwise from \eqref{desgamma} we have 
\[\frac{M}{|\g|^{as}-1} \geq \frac{1}{2} \ \Rightarrow \ s \leq \frac{\log(2M+1)}{a\log|\gamma|},\]
which implies the result. Back to \eqref{case2.2}, since $|\Gamma|<1/2$, then
\begin{equation*}
\frac{1}{2} < |chg^{-s}\g^{t-kas}| < \frac{3}{2} \ \Rightarrow \ |\log|ch| -s\log|g| + (t-kas)\log |\g|| < \log 2,
\end{equation*} 
and we get,
\[
|t-kas| < s \frac{\log|g|}{\log|\g|} + \frac{\log|2ch|}{\log|\g|}.
\]
Hence $|t-kas| = O(s)$, as we claimed. 

Picking up from \eqref{case2.2}, note that we can also \emph{look} at $\Gamma$ as 
\begin{equation}\label{GammaLinearForm}
    \Gamma = chg^{-s}\g^{t-kas}-1,
\end{equation}
with the aim to use Lemma \ref{lemma1} to get a lower bound for its absolute value. But before we are able to do that, we must consider the case where $t-aks$ equals zero. In order to do that, we have to split it into a few cases depending on $|g|$. First, let $|g|>1$. If $t=kas$, \eqref{GammaLinearForm} combined with $|\Gamma|<1/2$ gives,
\begin{equation} \label{Gamma2.2}
    |chg^{-s}- 1| < \frac{1}{2} \ \Rightarrow \ |chg^{-s}| \geq \frac{1}{2},
\end{equation}
then for $|g|>1$, we have the following upper bound for $s$ 
\[ 
s\le \frac{\log |2ch|}{\log |g|}, 
\]
thus we can take $E:= \log |2ch|/\log|g|$.
Now, we deal with $|g|<1$. Analogously, 
\begin{equation*}
    |chg^{-s}|\ge \frac{3}{2} < 2 \ \Rightarrow \ \log|ch| + s\log|g|^{-1} < \log 2,
\end{equation*} 
so,
 \[
 s\le \frac{\log2 - \log|ch|}{\log |g|^{-1}},
 \]
and we are done by taking $E:= (\log2 - \log|ch|)/\log |g|^{-1}$.
Now, if $|g| = 1$ we have two cases to deal with:
\begin{itemize}
    \item $\mathbf{|ch|\ne 1:}$ By \eqref{Gamma2.2}, we have 
\[|1-|ch|| \le \frac{M}{|\g|^{as}-1},\]
thus
\[\ s\le \frac{\log\left(1+\frac{M}{|1-|ch||}\right)}{a\log |\g|},\]
and we can take $E:= \log\left(1+\frac{M}{|1-|ch||}\right)/(a\log |\g|)$.
\item $\mathbf{|ch| = 1:}$ By \eqref{GammaLinearForm}, $|\Gamma +1|=|ch|=1$, and we recall that $\Gamma$ is a real number, so $|\Gamma|$ must be $0$ or $2$, but $|\Gamma|<1/2$, so $|\Gamma|=0$.  Now, since solving the case $|\Gamma| = 0$,  will be useful for us later to apply Lemma \ref{lemma1}, and requires a bit more work, we detail it next. First, notice that if we define
\[j_0 := \max\{0\leq j \leq k-1;\e_j \ne 0\},\]  
then,
\[|\Gamma|=0 \ \Rightarrow \ \e_0+\e_1\g^{as}+\cdots+\e_{j_0-1}\g^{as(j_0-1)}+\e_{j_0}\g^{asj_0} = 0,\]
By the triangle inequality and having in mind that $|\e_{j_0}|\geq m$, we get the following chain
\[
m|\g|^{asj_0}\le |\e_{j_0}||\g^{asj_0}| \le |\e_0|+|\e_1||\g|^{as}+\cdots+|\e_{j_0-1}||\g|^{as(j_0-1)},
\]
which gives,
\[m|\g|^{asj_0} \le  M(1+|\g|^{as}+\cdots+|\g|^{as(j_0-1)}) \leq  M \frac{|\g|^{asj_0}}{|\g|^{as}-1}.\]
Hence, taking logarithms, again we obtain an upper bound for $s$:
\[s\le \frac{\log \left(\frac{M+m}{m}\right)}{a\log|\g|}, \]
and by considering $E:= \log \left(\frac{M+m}{m}\right)/(a\log|\g|)$, we are done.
\end{itemize}

Finally, we are ready to apply Lemma \ref{lemma1}. Recall from \eqref{GammaLinearForm} that we must take
\[
\alpha_1:=ch \ ; \ \alpha_2:= g \ ; \ \alpha_3:= \g \ ; \]
\[
b_1:=1 \ ; \ b_2:= -s \ ; \ b_3:= t-kas.
\]
Note that $D = [\Q(\a_1,\a_2,\a_3):\Q] = [\Q(ch,g,\g):\Q]$, is a constant that depends only on $c$ and the parameters of $(G_n)_n$ and $(H_n)_n$, and the same holds for  
$A_j := D h(\alpha_j)+|\log(\alpha_j)|+0.16$, for $j=1,2,3$. Now, if we take $B:= s+|t-kas|$, then from our claim we have $B = O(s)$. We now compare our bounds on $|\Gamma|$ given by \eqref{desgamma} and Lemma \ref{lemma1}, to get

\[\frac{M}{|\g|^{as}-1}\ge |\Gamma| > \exp(-C_1A_1A_2A_3\log(eB)),\]
thus
\[|\g|^{as} < 1+M\exp(C_1A_1A_2A_3\log(eB)) < 2M\exp(C_1A_1A_2A_3\log(eB)),\]
and since $B = O(s)$, we have at last,
\[ as\log |\g|\le \log(2M)+C_1A_1A_2A_3\log(eB) = O(\log s),\]
which cannot hold for $s\gg 0$. Therefore, as stated, there is an effectively computable constant $E$, in this case depending only on the previously mentioned parameters, which finishes our proof. 
\qed



\begin{thebibliography}{9999}
\bibitem{gersica} Bedna$\check{\text{r}}\acute{\text{\i}}$k, G. Freitas, D. Marques and P. Trojovsk$\acute{\text{y}}$, On the sum of squares of consecutive $k$-bonacci numbers which are $l$-bonacci numbers, \emph{Colloq. Math.} \textbf{156} (2019), 153--164.

\bibitem{chaves} A.P. Chaves and D. Marques, A Diophantine equation related to the sum of squares of consecutive $k$-generalized Fibonacci numbers, \emph{Fibonacci Quart}. {\bf 52} (2014), no. 1, 70--74.

\bibitem{CMT} A.P. Chaves, D. Marques and A. Togb\' e, On the sum of powers of terms of a linear recurrence sequence. \emph{Bull. Braz. Math. Soc. (N.S.)}. \textbf{43} (2012), no. 3, 397--406.

\bibitem{kalman} D. Kalman and R. Mena, The Fibonacci numbers exposed, \textit{Math. Mag.} \textbf{76} (2003), no. 3, 167--181.


\bibitem{luca} F. Luca and R. Oyono, An exponential Diophantine equation related to powers of two consecutive Fibonacci numbers. \textit{Proc. Japan Acad. Ser}. A, \textbf{87} (2011) p. 45--50.

\bibitem{MT} D. Marques and A. Togb\' e, On the sum of powers of two consecutive Fibonacci numbers. \textit{Proc. Japan Acad. Ser. A,} \textbf{86} (2010) p. 174--176.


\bibitem{matveev} E. M. Matveev, An explict lower bound for a homogeneous rational linear form in logarithms of algebraic numbers, II, {\it Izv. Ross. Akad. Nauk Ser. Mat.} {\bf 64} (2000), 125--180. English transl. in {\it Izv. Math.} {\bf 64} (2000), 1217--1269.


\bibitem{pos} A. S. Posamentier, I. Lehmann, \textit{The (fabulous) Fibonacci numbers,} Prometheus Books, Amherst, NY, 2007.

\bibitem{ruiz}  C.A.G. Ruiz and F. Luca, An exponential Diophantine equation related to the sum of powers of two consecutive $k$-generalized Fibonacci numbers, \emph{Colloq. Math.} {\bf 137} (2014), 171--188.


\bibitem{shorey} T. N. Shorey and R. Tijdeman, \textit{Exponential Diophantine Equations}, Cambridge Tracts in Mathematics \textbf{87}, Cambridge University Press, Cambridge, 1986.

\end{thebibliography}
\end{document}